\documentclass[leqno, fleqn]{amsart}%
\usepackage{amsmath}
\usepackage{amsfonts}
\usepackage{amssymb}%

\makeatletter \DeclareMathSymbol{\Gamma}{\mathalpha}{letters}{"00}
\DeclareMathSymbol{\Delta}{\mathalpha}{letters}{"01}
\DeclareMathSymbol{\Theta}{\mathalpha}{letters}{"02}
\DeclareMathSymbol{\Lambda}{\mathalpha}{letters}{"03}
\DeclareMathSymbol{\Xi}{\mathalpha}{letters}{"04}
\DeclareMathSymbol{\Pi}{\mathalpha}{letters}{"05}
\DeclareMathSymbol{\Sigma}{\mathalpha}{letters}{"06}
\DeclareMathSymbol{\Upsilon}{\mathalpha}{letters}{"07}
\DeclareMathSymbol{\Phi}{\mathalpha}{letters}{"08}
\DeclareMathSymbol{\Psi}{\mathalpha}{letters}{"09}
\DeclareMathSymbol{\Omega}{\mathalpha}{letters}{"0A}
\DeclareMathSymbol{\varGamma}{\mathalpha}{operators}{"00}
\DeclareMathSymbol{\varDelta}{\mathalpha}{operators}{"01}
\DeclareMathSymbol{\varTheta}{\mathalpha}{operators}{"02}
\DeclareMathSymbol{\varLambda}{\mathalpha}{operators}{"03}
\DeclareMathSymbol{\varXi}{\mathalpha}{operators}{"04}
\DeclareMathSymbol{\varPi}{\mathalpha}{operators}{"05}
\DeclareMathSymbol{\varSigma}{\mathalpha}{operators}{"06}
\DeclareMathSymbol{\varUpsilon}{\mathalpha}{operators}{"07}
\DeclareMathSymbol{\varPhi}{\mathalpha}{operators}{"08}
\DeclareMathSymbol{\varPsi}{\mathalpha}{operators}{"09}
\DeclareMathSymbol{\varOmega}{\mathalpha}{operators}{"0A}

\newcommand{\allmodesymb}[2]{\relax\ifmmode{\mathchoice
{\mbox{\fontsize{\tf@size}{\tf@size}#1{#2}}}
{\mbox{\fontsize{\tf@size}{\tf@size}#1{#2}}}
{\mbox{\fontsize{\sf@size}{\sf@size}#1{#2}}}
{\mbox{\fontsize{\ssf@size}{\ssf@size}#1{#2}}}} \else
\mbox{#1{#2}}\fi}

\makeatother
\makeatletter
\renewcommand*\subjclass[2][2000]{%
  \def\@subjclass{#2}%
  \@ifundefined{subjclassname@#1}{%
    \ClassWarning{\@classname}{Unknown edition (#1) of Mathematics%
      Subject Classification; using '2000'.}%
  }{%
    \@xp\let\@xp\subjclassname\csname subjclassname@#1\endcsname%
  }%
} \makeatother

\theoremstyle{plain}

\theoremstyle{remark}

\allowdisplaybreaks \numberwithin{equation}{section}
\begin{document}
\title{Some Topics Related to Bergman Kernel}

\author{YIN Weiping }
\address{YIN Weiping: Dept. of Math., Capital Normal Univ., Beijing 100037, China}
\email{wyin@mail.cnu.edu.cn; wpyin@263.net}

\maketitle

Actually we will discuss some topics related to Bergman kernel on
Cartan-Hartogs domain. Cartan-Hartogs domain is introduced by Guy
Roos and Weiping YIN in 1998, which is built on the Cartan
domains(classical domains). The four big types of Cartan domains can
be written as[1]:$$
\begin{array}{ll}
R_I(m,n)&:=\{Z\in {\bf{C}^{mn}}: I-Z\overline{Z}^t>0, \},\\[3mm] R_{II}(p)&:=\{Z\in
{\bf{C}^{p(p+1)/2}}: I-Z\overline{Z}^t>0, \},\\[3mm] R_{III}(q)&:=\{Z\in
{\bf{C}^{q(q-1)/2}}: I-Z\overline{Z}^t>0,\},\\[3mm] R_{IV}(n)&:=\{Z\in {\bf{C}^n}:
1+|ZZ^t|^2-2Z\overline{Z}^t>0,\\& 1-|ZZ^t|^2>0\}.
\end{array}
$$
where $Z$ is $m\times n$ matrix, $p$ degree symmetric matrix, $q$
degree skew symmetric matrix and $n$ dimensional complex vector
respectively. Then the Cartan-Hartogs domains can be introduced as
follows:
$$
Y_I:=Y_I(N,m,n;K):=\{W\in {\bf{C}^N},Z\in R_I(m,n):$$
$$|W|^{2K}<\det(I-Z\overline{Z}^t),K>0\},
$$
$$Y_{II}:= Y_{II}(N,p;K):=\{W\in {\bf{C}^N},Z\in
R_{II}(p):$$ $$|W|^{2K}<\det(I-Z\overline{Z}^t),K>0\}, $$
$$
Y_{III}:=Y_{III}(N,q;K):=\{W\in {\bf{C}^N},Z\in R_{III}(q):$$
$$|W|^{2K}<\det(I-Z\overline{Z}^t),K>0\}, $$
$$Y_{IV}:=Y_{IV}(N,n;K):=\{W\in {\bf{C}^N},Z\in
R_{IV}(n):$$ $$|W|^{2K}<1-2Z\overline{Z}^t+|ZZ^t|^2, K>0\},
$$
where  $\det$ indicates the determinant, $N, m, n, p, q$ are natural
numbers. These domains are also called super-Cartan domains.

If the right hand of above inequalities are denoted by the
$N_j:=N_j(Z,\overline Z), j=I,II,III,IV$ respectively, then the
definition of Cartan-Hartogs domain can be also written as
$$
Y_j=\{W\in {\bf{C}^N}, Z\in R_j: |W|^{2K}< N_j(Z,\overline Z) \},
j=I,II,III,IV.$$ The following topics will be discussed:

I. The zeroes of Bergman kernel of Cartan-Hartogs domain;

II. The classical (Cannonical) metrics (Bergman metric, Caratheodory
metric, Kaehler-Einstein metric, Kobayashi metric) on Cartan-Hartogs
domain, which contains Bergman metric equivalent to Kaehler-Einstein
metric, Lu Qikeng constant, Bergman Kaehler-Einstein metric and some
good new metrics.

III. Generalized Cartan-Hartogs domain;

IV. The centre of representative domain and applications;

V. The solution of Dirichlet's problem of complex Monge-Amp\`ere
equation on Cartan-Hartogs domain and Kaehler-Einstein metric with
explicit formula.

\section*{I. The zeroes of Bergman kernel on Cartan-Hartogs domain}

The Cartan-Hartogs domain of the first type is defined by
$$
Y_I(N,m,n;K)=\{W\in {\bf{C^N}},Z\in
R_I(m,n):|W|^{2K}<det(I-Z\overline{Z}^t),K>0\}.$$ And
$$Y_I(1,1,n;K)=\{W\in {\bf{C}}, Z\in {\bf{C^n}}: |W|^{2K}+|z_1|^2+|z_2|^2+\dots+|z_n|^2<1\}.$$
Then the Bergman kernel of $Y_I(1,1,n;K)$ is
$$
K_I((W,Z);(\zeta,\xi))=K^{-n}\pi^{-(n+1)}F(Y)
det(I-Z\overline{\xi}^t)^{-(1+n+1/K)}.
$$
Where $$ F(Y)=\sum_{i=0}^{n+1}b_i\Gamma(i+1)Y^{i+1}, Y=(1-X)^{-1},
X=W\zeta[det(I-Z\overline{\xi}^t)]^{-1/K}.
$$ And $b_0=0$, let
$$
P(x)= (x+1)[(x+1+Kn)(x+1+K(n-1))\dots(x+1+K)].
$$ Then the others
$b_i(i=1,2,\dots\dots,n+1)$ can be got by
$$
b_i=[P(-i-1)-\sum_{k=0}^{i-1}b_k(-1)^k\Gamma(i+1)/\Gamma(i-k+1)][(-1)^i\Gamma(i+1)]^{-1}.
\eqno{(1.1)}
$$
Recently, Liyou Zhang prove that above formula can be rewritten as
$$b_i=\sum_{j=1}^{i}\frac{(-1)^jP(-j-1)}{\Gamma(j+1)\Gamma(i-j+1)}.$$

It is well known that for the first type of Cartan-Hartogs domains
there exists the holomorphic automorphism $(W^*, Z^*)=F(W,Z)$ such
that $F(W,Z_0)=(W^*,0)$ if $(W,Z_0)\in Y_I$. Due to the
transformation rule of Bergman kernel, one has $$
K_I((W,Z);(\zeta,\xi))=(detJ_F(W,Z))|_{Z_0=Z}K_I[(W^*,0),(\zeta^*,\xi^*)](det\overline{J_F(\zeta,\xi)}).$$
Therefore the zeroes of $K_I((W,Z);(\zeta,\xi))$ are same as the
zeroes of $K_I[(W^*,0),(\zeta^*,\xi^*)]=K^{-n}\pi^{-(n+1)}F(Y)$. Let
$W^*$ be the $W$, and $\zeta^*$ be the $\zeta$, then we have
$$K^{-n}\pi^{-(n+1)}F(Y)=K^{-n}\pi^{-(n+1)}F(y),
y=(1-W\overline{\zeta})^{-1}.\eqno{(1.2)}$$ Where
$$F(y)=\sum_{i=0}^{n+1}b_i\Gamma(i+1)y^{i+1}. \eqno{(1.3)}
$$ If $(W,0)$, $(\zeta, 0)$, $(W^*,0)$ and $(\zeta^*, 0)$  belong to
$Y_I$, then their norms $ |W|,|\zeta|,|W^*|,|\zeta^*|$ are less than
1.

1.1. Let $t=W\overline{\zeta}$, then $|t|<1$, and
$F(y)=(1-t)^{-(n+2)}G(t)$, where
$G(t)=\sum_{i=0}^{(n+1)}b_i\Gamma(i+1)(1-t)^{n+1-i}$. {\bf Therefore
to discuss the the presence or absence of zeroes of the Bergman
kernel function of $Y_I(1,1,n;K)$ can be reduced to discuss the
zeroes of polynomial with real coefficients in the unit disk in
$\bf{C}$}[3].

Because $y=(1-W\overline{\zeta})^{-1}=\frac{1}{1-t}$, which maps the
unit disk in $t$-plane onto the half-plane in $y$-plane $Re y>1/2$,
{\bf therefore to discuss the the presence or absence of zeroes of
the Bergman kernel function of $Y_I(1,1,n;K)$ can be reduced to
discuss the zeroes of polynomial with real coefficients in the right
half-plane $Re y>1/2$.}

Above two statements are true not only for the $Y_I(1,1,n;K)$ but
also for all of the Cartan-Hartogs domains(and Hua domains).

1.2. In the low dimension case, it is very easy to answer the Lu
Qikeng problem for the Cartan-Hartogs domain. For example, we can
say that $Y_I(1,1,1;K)$ is Lu Qi-Keng domain.

At this time $Y_I(1,1,1;K)=\{W\in {\bf{C}}, Z\in {\bf{C}}:
|W|^{2K}+|Z|^2<1\}$, and the zeroes of its Bergman kernel function
$K_I[(W,Z),(\zeta,\xi)]$ are same as the zeroes of
$K_I[(W^*,0),(\zeta^*,0)]$. But
$$K_I[(W,0),(\zeta,0)]=K^{-n}\pi^{-(2)}F(y), F(y)=\sum_{i=0}^{2}b_i\Gamma(i+1)y^{i+1}$$
$$y=(1-t)^{-1}.$$
Where $b_1=K-1, b_2=1, b_0=0$, therefore
$$F(y)=(K-1)y^2+2y^3=y^3[(K-1)(1-t)+2].$$
But the zeroes of $F(y)$ are equal to $t=(K+1)/(K-1)$, its norm
$|t|>1$, it is impossible.

{\bf Therefore the Bergman kernel function of $Y_I(1,1,1;K)$ is
zero-free, that is the $Y_I(1,1,1;K)$ is Lu Qi-Keng domain.}
Therefore we also prove that:

{\bf If $D\subset{C^2}$ is a bounded pseudoconvex domain with real
analytic boundary and its holomorphic automorphism group is
noncompact, then $D$ is the Lu Qi-Keng domain} due to the E.Bedford
and S.I.Pinchuk's following theorem[4].

{\bf Theorem 1.1:} If $D$ is a bounded pseudoconvex domain with real
analytic boundary and its holomorphic automorphism group is
noncompact, then $D$ is biholomorphically equivalent to a domain of
the form$$ E_m=\{(z_1,z_2)\in{\bf{C^2}}: |z_1|^{2m}+|z_2|^2<1\}$$
for some positive integer $m$.

\section*{II. The classical (Cannonical) metrics on Cartan-Hartogs domain}

Let $D$ be the bounded domain in $\bf{C^M}$, $B_D, C_D, KE_D, K_D$
denote the Bergman metric, Caratheodory metric, Kaehler-Einstein
metric, Kobayashi metric respectively, then we have $C_D\leq B_D,
C_D\leq K_D$ [5], and if $D$ is also the convex domain then
$C_D=K_D$[5]. On the other hand, there is no clear relationship
between the $B_D$ and $K_D$.

2.1. {\bf But we have that $B_D\leq cK_D$ for the Cartan-Hartogs
domain where $c$ is the constant} [6-9].

2.2. {\bf We proved the Bergman metric is equivalent to
Kaehler-Einstein metric[10], that is $B_D\sim KE_D$ on
Cartan-Hartogs domains}. For example, we consider the Cartan-Hartogs
of the first type $ Y_I=Y_I(N,m,n;K)$. Let $$
G_{\lambda}=G_{\lambda}(Z,W)=Y^{\lambda}[\det(I-Z\overline
Z^t)]^{-(m+n+\frac{N}{K})}, \lambda > 0,
$$
$$T_{\lambda I}(Z,W;\overline{Z},\overline{W})
=(g_{i\overline{j}})=\left(\frac{\partial^2\log G_{\lambda}}
{\partial z_i\partial \overline{z}_{j}}\right),$$ where
$$Y=(1-X)^{-1}, X=|W|^2[\det(I-Z\overline Z^t)]^{-\frac 1K}.$$
Then $G_{\lambda}$ induces a metric
$$Y(I\lambda):=[(dw,dz)T_{\lambda I}(Z,W;\overline{Z},\overline{W})\overline{(
dw,dz)}^t]^{1/2}.$$

{\bf Firstly, by the direct computations one can prove that
$B_{Y_I}\sim Y(I\lambda)$. The $Y(I\lambda)$ has good properties:
Its holomorphic sectional curvature and Ricci curvature are bounded
from above and below by the Negative constants. Then based on above
good properties and using the Yau's Schwarz lemma[11] one can prove
$KE_{Y_I}\sim Y(I\lambda)$. Therefore $KE_{Y_I}\sim B_{Y_I}$, and
the metric $Y(I\lambda)$ may be useful for us.}

2.3. {\bf Definition:} A complex manifold $M^n$ is called
holomorphic homogeneous regular if there are positive constants
$r<R$ such that for each point $p\in M$, there is a one to one
holomorphic map $f: M\longrightarrow \bf{C}^n$ such that

i) $f(p)=0;$

ii) $B_r\subset f(M)\subset B_R,$ where $B_r$ and $B_R$ are balls
with radius $r$ and $R$ respectively.

 {\bf Theorem 2.1(Liu-Sun-Yau)}[11,12]: For holomorphic homogeneous regular manifolds,
the Bergman metric, the Kobayashi metric and the Caratheodory metric
are equivalent.

{\bf Therefore if Cartan-Hartogs domains are the holomorphic
homogeneous regular manifolds, then the Bergman metric, the
Kobayashi metric, the Caratheodory metric and Kaehler-Einstein
metric are equivalent. But whether the Cartan-Hartogs domains are
the holomorphic homogeneous regular manifolds? This problem remains
open.}

2.4. From an immediate consequence of an inequality due to Lu
Qikeng's paper [13], we have the following

{\bf Theorem 2.2(Lu Qikeng):} Let $D$ be a bounded domain in
$\bf{C}^n$. Then for each tangent $v\in T_z(D)=\bf{C}^n$ at $z\in
D$, $B_D(z,v)\geq C_D(z,v)$. Where $B_D(z,v)$ equals the length of
$v$ w.r.t. the Bergman metric $B_D$, and $C_D(z,v)$ equals the
length of $v$ w.r.t. the differential Caratheodory metric $C_D$.

{\bf Therefore Cheung and Wong introduce the definition of Lu
constant $L(D)$ of a bounded domain $D$ in $\bf{C}^n$ as follows[5].

Definition:
$$L(D)= \sup\limits^{z\in D}_{v\not=0\in{T_z(D)}}
(\frac{C_D(z,v)}{B_D(z,v)})$$}

Lu's theorem says that $L(D)\leq 1$. $L(D)=(1/(n+1))^{1/2}$ when $D$
is the unit ball in $\bf{C}^n$. One can try to determine the Lu's
constants of all Cartan domains and all Cartan-Hartogs domains.

2.5. {\bf Some years ago S.T.Yau proposed an intricate problem to
look for a characterization of the bounded pseudoconvex domains on
which the Bergman metrics are complete Kaehler-Einstein metric[5].}

The following Lu's theorem can be viewed  as a particular case of
Yau's problem of which the Bergman metric is of constant negative
holomorphic sectional curvature:

{\bf Theorem 2.3(Lu Qikeng)}[14]: Let $D$ be a bounded domain in
$\bf{C}^n$ with a complete Bergman metric $B_D$. If the holomorphic
sectional curvature is equal to a negative constant $-c^2$, then $D$
is biholomorphic to the Euclidean ball $B_n$ and
$c^2=\frac{2}{n+1}$.

S.Y.Cheng conjectures that {\bf a strangle pseudoconvex domain whose
Bergman metric is Kaehler-Einstein must be biholomorphic to the
Euclidean ball}[5].

{\bf We can also prove that if the Bergman metric of Cartan-Hartogs
domain is Kaehler-Einstein, then this Cartan-Hartogs domain must be
homogeneous(See below).}

\section*{III. Generalized Cartan-Hartogs domain}

Some years ago we generate the Cartan-Hartogs domain to the Hua
domain as follows[3]:
$$
\begin{array}{lll}
\{ W_j\in {\bf{C}^{N_j}},Z\in
{R_s}:\displaystyle\sum_{j=1}^r\frac{||W_j||^{2p_j}}{[N_s(Z,\overline
Z)]^{K_j}}<1, \\p_j>0, K_j>0, j=1,\dots,r \}.
s=I,II,III,IV.\end{array}$$

{\bf Right now we will generate the Cartan-Hartogs domain from
another way.}

Let $\Omega$ be a domain in $\bf{C}^n$, $\rho$ a positive continuous
function on $\Omega$, and let $D$ be a (fixed) irreducible bounded
symmetric domain in $\bf{C}^d$. Then Roos, Engli$\breve{s}$ and
Zhang define a new domain in $\bf{C}^{n+d}$ as follows[15,16]:
$$
\Omega^{D}:=\{(w,z)\in {\bf{C}^d \times \Omega}:
\frac{w}{\rho(z)}\in D \}.\eqno{(3.1)}$$

Let $B(0,1)$ be the unit ball in $\bf{C}^d$, and let
$$D=B(0,1), \rho(z)=N_j(z,z)^{1/(2K)},$$ then one has
$$\Omega^{B(0,1)}:=\{(w,z)\in {\bf{C}^d \times \Omega }:
\frac{w}{N_j(z,z)^{1/(2K)}}\in B(0,1) \}.\eqno{(3.2)}$$

The $\frac{w}{N_j(z,z)^{1/(2K)}}\in B(0,1)$ can be denoted by
$$(\frac{w}{N_j(z,z)^{1/(2K)}})\overline{(\frac{w}{N_j(z,z)^{1/(2K)}})}'<1.$$
That is
$$|w|^{2K}<N_j(z,z).$$Therefore
$$\Omega^{B(0,1)}=\{(w,z)\in {\bf{C}^d \times \Omega }:
|w|^{2K}<N_j(z,z) \}. \eqno{(3.3)}$$ Above (3.3) is the definition
of Cartan-Hartogs domain.

Let
$$D=R_{j}, \rho(z)=N_i(z,z)^{1/(2K)}, d=dim D, \Omega=R_i,$$
then we get the following new domain, which generalizes the
Cartan-Hartogs, and is called generalized Cartan-Hartogs domain:
$$R_i^{R_j}=\{(w,z)\in {\bf{C}^d \times R_i}:
\frac{w}{N_i(z,z)^{1/(2K)}}\in R_j \}. \eqno{(3.4)}$$ where
$i,j=I,II,III,IV$. Therefore we get 16 types of generalized
Cartan-Hartogs domain as follows: $$Y(I,I)=R_I^{R_I}=\{(w,z)\in
{\bf{C}^{mn}\times R_I }:
w\overline{w'}<det(I-z\overline{z'})^{1/(K)}I^{(m)}\}.$$
$$Y(I,II)=R_{I}^{R_{II}}=\{(w,z)\in {\bf{C}^{p(p+1)/2}\times R_{I}
}: w\overline{w'}<det(I-z\overline{z'})^{1/(K)}I^{(p(p+1)/2)}\}.$$
$$Y(I,III):=R_{I}^{R_{III}}=\{(w,z)\in {\bf{C}^{q(q-1)/2}\times R_{I}  }:
w\overline{w'}<det(I-z\overline{z'})^{1/(K)}I^{(q(q-1)/2)}\}.$$
$$Y(I,IV):=R_{I}^{R_{IV}}=\{(w,z)\in {\bf{C}^{N}\times R_{I}  }:
2det(I-z\overline{z'})^{1/(K)}w\overline{w'}-|ww'|^2<det(I-z\overline{z'})^{2/(K)},$$
$$|ww'|<det(I-z\overline{z'})^{1/(K)}\}.$$
$$Y(II,I):=R_{II}^{R_{I}}=\{(w,z)\in {\bf{C}^{mn}\times R_{II}  }:
w\overline{w'}<det(I-z\overline{z'})^{1/(K)}I^{(m)}\}.$$
$$Y(II,II):=R_{II}^{R_{II}}=\{(w,z)\in {\bf{C}^{p(p+1)/2}\times R_{II}  }:
w\overline{w'}<det(I-z\overline{z'})^{1/(K)}I^{(p)}\}.$$
$$Y(II,III):=R_{II}^{R_{III}}=\{(w,z)\in {\bf{C}^{q(q-1)/2}\times R_{II}  }:
w\overline{w'}<det(I-z\overline{z'})^{1/(K)}I^{(q)}\}.$$
$$Y(II,IV):=R_{II}^{R_{IV}}=\{(w,z)\in {\bf{C}^{N}\times R_{II}  }:
2det(I-z\overline{z'})^{1/(K)}w\overline{w'}-|ww'|^2<det(I-z\overline{z'})^{2/(K)},$$
$$|ww'|<det(I-z\overline{z'})^{1/(K)}\}.$$
$$Y(III,I):=R_{III}^{R_{I}}=\{(w,z)\in {\bf{C}^{mn}\times R_{III}  }:
w\overline{w'}<det(I-z\overline{z'})^{1/(K)}I^{(m)}\}.$$
$$Y(III,II):=R_{III}^{R_{II}}=\{(w,z)\in {\bf{C}^{p(p+1)/2}\times R_{III}  }:
w\overline{w'}<det(I-z\overline{z'})^{1/(K)}I^{(p)}\}.$$
$$Y(III,III):=R_{III}^{R_{III}}=\{(w,z)\in {\bf{C}^{q(q-1)/2}\times R_{III}  }:
w\overline{w'}<det(I-z\overline{z'})^{1/(K)}I^{(q)}\}.$$
$$Y(III,IV):=R_{III}^{R_{IV}}=\{(w,z)\in {\bf{C}^{N}\times R_{III}  }:
2det(I-z\overline{z'})^{1/(K)}w\overline{w'}-|ww'|^2<det(I-z\overline{z'})^{2/(K)},$$
$$|ww'|<det(I-z\overline{z'})^{1/(K)}\}.$$
$$Y(IV,I):=R_{IV}^{R_{I}}=\{(w,z)\in {\bf{C}^{mn}\times R_{IV}  }:
w\overline{w'}<(1+|zz'|^2-2z\overline{z'})^{1/(K)}I^{(m)}\}.$$
$$Y(IV,II):=R_{IV}^{R_{II}}=\{(w,z)\in {\bf{C}^{p(p+1)/2}\times R_{IV}  }:
w\overline{w'}<(1+|zz'|^2-2z\overline{z'})^{1/(K)}I^{(p)}\}.$$
$$Y(IV,III):=R_{IV}^{R_{III}}=\{(w,z)\in {\bf{C}^{q(q-1)/2}\times R_{IV}  }:
w\overline{w'}<(1+|zz'|^2-2z\overline{z'})^{1/(K)}I^{(q)}\}.$$
$$Y(IV,IV):=R_{IV}^{R_{IV}}=\{(w,z)\in {\bf{C}^{N}\times R_{IV}  }:
2(1+|zz'|^2-2z\overline{z'})^{1/(K)}w\overline{w'}-|ww'|^2<(1+|zz'|^2-2z\overline{z'})^{2/(K)},$$
$$|ww'|<(1+|zz'|^2-2z\overline{z'})^{1/(K)}\}.$$
{\bf These are the new research fields, one can seeks the Bergman
kernel, Szeg$\ddot{o}$ and consider other topics.}

\section*{IV. The centre of representative domain and applications }

The Riemann mapping theorem characterizes the planar domains that
are biholomorphically equivalent to the unit disk. In the higher
dimensions, there is no Riemann mapping theorem, and the following
problem arise:

Are there canonical representatives of biholomorphic equivalence
classes of domains?

In the dimension one, if $K(z,w)$ is the Bergman kernel function of
simply connected domain $D\neq {\bf{C}}$, it is well known that the
biholomorphic mapping
$$F(z)=\frac{1}{K(t,t)}\frac{\partial}{\partial \overline{w}}log\frac{K(z,w)}
{K(w,w)}|_{w=t}$$ maps the $D$ onto unit disk.

In the higher dimensions, Stefan Bergman introduced the notion of a
"representative domain" to which a given domain may be mapped by
"representative coordinates". If $D$ is a bounded domain in
${\bf{C^n}}$, $K(Z,W)$ is the Bergman kernel function of $D$, let
$$T(Z,Z)=(g_{ij})=(\frac{\partial^2logK(Z,W)}{\partial z_i\partial \overline{z_j}})$$ and its converse is
$T^{-1}(Z,W)=(g^{-1}_{ji})$. Then the local representative
coordinates based at the point $t$ is
$$f_i(Z)=\sum^n_{j=1}g^{-1}_{ji}\frac{\partial}{\partial{\overline{W_j}}}
log\frac{K(Z,W)}{K(W,W)}\mid_{W=t}, i=1,\dots,n.$$ Or

$$F(Z)=(f_1,\dots,f_n)=\frac{\partial}{\partial{\overline{W}}}
log\frac{K(Z,W)}{K(W,W)}\mid_{W=t}T^{-1}(t,t),$$ where
$$\frac{\partial}{\partial{\overline{W}}}=(\frac{\partial}{\partial{\overline{W_1}}},
\frac{\partial}{\partial{\overline{W_2}}}, \dots,
\frac{\partial}{\partial{\overline{W_n}}}).$$ These coordinates take
$t$ to $0$ and have complex Jacobian matrix at $t$ is equal to the
identity. {\bf The $F(D)$ is called the representative domain of
$D$. If $D$ is biholomorphic equivalent to $D_1$, then $D$ and $D_1$
have same representative domain.}

Zeroes of the Bergman kernel function $K(Z,W)$ evidently pose an
obstruction to the global definition of Bergman representative
coordinates. This observation was Lu Qi-Keng's motivation for asking
which domains have zero-free Bergman kernel functions. This problem
is called Lu Qi-Keng conjecture by M.Skwarczynski in 1969 in his
paper [17]. {\bf If the Bergman kernel function of $D$ is zero-free,
that means the Lu Qi-Keng conjecture has a positive answer, then the
domain $D$ is called the Lu Qi-Keng domain.}

4.1. In 1981 Lu Qikeng introduces an another definition of
"representative domain"[18].

{\bf Definition: A bounded domain in $\bf{C}^n$ is called a
representative domain, if there is a point $t\in{D}$ such that the
matrix of the Bergman metric tensor $T(z,\overline{t})$ is
independent of $z\in{D}$. The point $t$ is called the centre of the
representative domain.}

If $D$ is representative domain in the sense of Lu, and $D_1$ is the
representative domain of $D$ in the sense of Bergman, then $D$ is
same as the $D_1$ under an affine transformation.

4.2. In 1981, Lu Qikeng[18] proved the following

{\bf Theorem 4.1:} Let $D$ be a bounded domain and $D_1$ be a
representative domain of Lu in $\bf{C}^n$ with centre $s_0$. If $f:
D\longrightarrow D_1$ is a biholomorphic mapping, then $f$ is of the
form
$$
f(z)=s_0+[\frac{\partial}{\partial
\overline{t}}log\frac{K(z,\overline{t})}
{K(t,\overline{t})}]_{t=t_0}T^{-1}(t_0,\overline{t_0})A.$$ Moreover
$K(z,\overline{t_0})$ is zero free when $z\in D$. Where
$s_0=f(t_0)$, $A=(\frac{\partial f}{\partial z})_{z=t_0}$, and
$\frac{\partial}{\partial \overline{t}}=(\frac{\partial}{\partial
\overline{t_1}},\dots,\frac{\partial}{\partial \overline{t_n}})$.

{\bf Corollary 1:} Let $D$ be a bounded domain and $D_1$ be a
representative domain of Lu in $\bf{C}^n$ with centre $s_0$, if $D$
is biholomorphic equivalent $D_1$, the Bergman kernel of $D_1$ is
$K_1(w,\overline{s})$, then $K_1(w,\overline{s_0})$ is zero free
when $w\in{D_1}$.

{\bf Corollary 2:} If $s_0=0$, $D=D_1$, then $A=I$, and the
holomorphic automorphism $f(z)$ of $D$ has the following form:
$$
f(z)=[\frac{\partial}{\partial
\overline{t}}log\frac{K(z,\overline{t})}
{K(t,\overline{t})}]_{t=t_0}T^{-1}(t_0,\overline{t_0}).$$ Where
$0=f(t_0)$.

4.3. If the full  group of holomorphic automorphism is denoted by
$Aut(D)$, and let $S=\{z: f(z)=0, f\in Aut(D) \}$. Then the $Aut(D)$
is constituted by
$$
f(z)=[\frac{\partial}{\partial
\overline{t}}log\frac{K(z,\overline{t})}
{K(t,\overline{t})}]_{t=t_0}T^{-1}(t_0,\overline{t_0}).$$ Where
$t_0$ spreads all over $S$.

{\bf Therefore if $0$ is the centre of representative domain $D$,
and the set $S$ is got explicitly, then the $Aut(D)$ can be got
explicitly as above. From this, we can get the $Aut(D)$ with
explicit form if $D$ is the Cartan-Hartogs domain.}

\section*{V. The solution of Dirichlet's problem of complex Monge-Amp\`ere
equation on Cartan-Hartogs domain and Kaehler-Einstein metric with
explicit formula}

Complex Monge-Amp\`ere equation is the nonlinear equation with high
degree, therefore to get its solution is very difficult.

S.Y.Cheng, N.M. Mok, S.T. Yau consider the following Dirichlet's
problem of the complex Monge-Amp\`ere equation:
$$
 \left\{
\begin{array}{ll}
\displaystyle {\rm{det}} \left( \frac{\partial^2g}{\partial z_i
\partial \overline{z}_j} \right) =e^{(n+1)g} & z \in  D, \\
\displaystyle g=\infty & z \in \partial D,\\
\end{array}
\right. $$ And they proved that the above problem exists unique
solution[19,20], where the $g$ can induce the Kaehler-Einstein
metric as follows: $$ (KE_D)^2=dz(\frac{\partial^2g}{\partial z_i
\partial \overline{z}_j} )\overline{dz}^t.$$

We consider the explicit solution of Dirichlet's problem of complex
Monge-Amp\`ere equation on $Y_I$:
$$
 \left\{
\begin{array}{ll}
\displaystyle {\rm{det}} \left( \frac{\partial^2g}{\partial z_i
\partial \overline{z}_j} \right)_{1\leq i,j\leq M} =e^{(M+1)g} & z \in  Y_I, \\
\displaystyle g=\infty & z \in \partial Y_I,\\
\end{array}
\right. \eqno{(5.1)}$$  where $M=N+mn$ is the complex dimension of
$Y_{I}$.

Because $Y_I$ is pseudoconvex domain. Therefore the solution of
problem (5.1) is existent and unique.

5.1. We prove that the solution of problem (5.1) can be got in
semi-explicit formula, and the explicit solution is obtained in
special case. That is the following theorem is proved[21]:

{\bf Theorem 5.1:} If $G(X)$ is the solution of the following
problem
$$\left\{
\begin{array}{lll}(M+1)^{-M}[\frac{X}{K}G'+(m+n+\frac{N}{K})G]^{mn}
[GG'+(GG''-(G')^2)X]\frac{(G')^{N-1}}{G^{M+1}}=G,\\
G(0)=\displaystyle K^{-mn}; \displaystyle lim_{X\rightarrow
1}G(X)=\infty,
\end{array} \right.\eqno{(5.2)}$$ then
$$g=(M+1)^{-1}log[G(X)det(I-Z\overline{Z}^t)^{-(m+n+N/K)}]$$ is
the solution of the problem (5.1); if $K=\frac{mn+1}{m+n}$, and
$$G(X)=(\frac{m+n}{mn+1})^{mn}(1-X)^{-(M+1)},$$ then the following
$g$ is the special solution of the problem (5.1):
$$g=(M+1)^{-1}log[(\frac{m+n}{mn+1})^{mn}(1-X)^{-(M+1)}det(I-Z\overline{Z}^t)^{-(m+n+N/K)}]$$
$$=log[(1-X)^{-1}det(I-Z\overline{Z}^t)^{-\frac{m+n}{mn+1}}(\frac{m+n}{mn+1})^{\frac{mn}{M+1}}],\eqno{(5.3)}$$
where
$$ X=X(Z,W)=|W|^2[det(I-Z\overline{Z}^t)]^{-1/K}, G'=\frac{d
G(X)}{d X}, G''=\frac{d^2G(X)}{d X^2}.$$

 {\bf Remark 1:} The complex Monge-Amp\`ere equation is the nonlinear
equation, hence to get its explicit solution is very difficult.
Therefore mathematicians hope to get the solutions for the problem
(5.1) by using the numerical method. Due to the above results the
numerical method of the problem (5.1) is reduced to the numerical
method of the problem (5.2). Which reduce the complexity of the
numerical method of problem (5.1) consumedly. Next, if the numerical
method of problem (5.1) or the numerical method of problem (5.2) is
appeared in the future, then the special solution $g$ (see (5.3))
can be used to check these numerical methods. And if one reduces the
complex Monge-Amp\`ere equation in (5.1) by the linearization
method, then the $g$ of (5.3) can be also to check the precision and
the rationality for the linearization method.

{\bf Remark 2:} Although the problem (5.1) have not been got the
explicit solution in general case, but its semi-explicit solution
has the form
$$g=(M+1)^{-1}log[G(X)det(I-Z\overline{Z}^t)^{-(mn+N/K)}],
$$   where $G(X)$ satisfies the (5.2).

{\bf Remark 3:}  The Bergman kernel function of $Y_I(N,m,n;K)$ is
$$
K_I(W,Z;\overline{W},\overline{Z})=K^{-mn}\pi^{-(mn+N)}G(X)
det(I-Z\overline{Z}^t)^{-(m+n+N/K)}.
$$
Where
$$
G(X)=\sum_{i=0}^{mn+1} b_i\Gamma(N+i)(1-X)^{-(N+i)},
X=X(W,Z)=|W|^2[det(I-Z\overline{Z}^t)]^{-1/K},$$
$$
|W|^2=\sum_{j=1}^N|W_j|^2,
$$
and $b_i$ are constants.

{\bf If $Y_I$ is homogeneous domain, then its Bergman metric is
equal to its Kaehler-Einstein metric, that is
$G(X)=\sum_{i=0}^{mn+1} b_i\Gamma(N+i)(1-X)^{-(N+i)}$ must satisfies
the equation (5.2). If $G(X)$ is not satisfies (5.2), then $Y_I$ is
not homogeneous.}

{\bf If $G(X)$ is satisfies (5.2), then the Bergman metric of $Y_I$
is equal to the Kaehler-Einstein metric.}

By computations, {\bf we prove that the $G(X)$ satisfies the
equation (5.2) if and only if $m=1$.That is the $Y_I$ is the unit
ball(homogeneous domain).}

\end{document}